\documentclass[12pt]{article}
\usepackage{amsthm,amsmath,amssymb}

\usepackage{graphicx}
\usepackage[T1]{fontenc}

\usepackage{algorithm2e}
\usepackage[colorlinks=true,citecolor=black,linkcolor=black,urlcolor=blue]{hyperref}


\newcounter{parentnumber}
\theoremstyle{plain}
\newtheorem{theorem}{Theorem}

\theoremstyle{lemma}
\newtheorem{lemma}{Lemma}

\theoremstyle{proposition}
\newtheorem{proposition}{Proposition}

\theoremstyle{corollary}
\newtheorem{corollary}{Corollary}

\theoremstyle{definition}
\newtheorem{definition}{Definition}
\newtheorem{example}[theorem]{Example}
\newtheorem{conjecture}[theorem]{Conjecture}

\theoremstyle{remark}
\newtheorem{remark}[theorem]{Remark}

 \title{\bf  A Simple Explanation for the Reconstruction of Graphs\\ 
 \small{Second version}  }

\author{Ameneh Farhadian\\ %\thanks{The contents of this paper are taken from the author's Ph.D. Thesis, Department of Mathematical Sciences, Sharif University of Technology, Supervised by Professor E. S. Mahmoodian} \\
\small Department of Mathematical Sciences\\[-0.8ex]
\small Sharif University of Technology\\[-0.8ex]
 }

\date{ \small Mathematics Subject Classifications:  05C60}

\begin{document}

\maketitle

% E-JC papers must include an abstract. The abstract should consist of a
% succinct statement of background followed by a listing of the
% principal new results that are to be found in the paper. The abstract
% should be informative, clear, and as complete as possible. Phrases
% like "we investigate..." or "we study..." should be kept to a minimum
% in favor of "we prove that..."  or "we show that...".  Do not
% include equation numbers, unexpanded citations (such as "[23]"), or
% any other references to things in the paper that are not defined in
% the abstract. The abstract will be distributed without the rest of the
% paper so it must be entirely self-contained.

\begin{abstract}
The graph reconstruction conjecture states that all graphs on at least three vertices are determined up to isomorphism by their deck.
In this paper, a general framework for this problem is proposed to simply explain the reconstruction of graphs. Here, we do not prove or reject the reconstruction conjecture. But, we explain why a graph is reconstructible. For instance, the reconstruction of small graphs which have been shown  by computer, is explained in this framework.\\
We show that any non-regular graph  has a proper induced subgraph which is  unique due to either its structure or the way of its connection to  the rest of the graph. Here, the former subgraph is defined  an anchor and the latter a connectional anchor, if it is distinguishable in the deck.
We show that if a graph has an orbit with at least three vertices whose removal leaves an anchor, or  it has two vertices whose removal leaves an anchor with the mentioned condition in the paper, then it is reconstructible. This  simple statement can easily explain the reconstruction of a graph from its deck.

\end{abstract}

\section {Introduction}
The Reconstruction Conjecture is an interesting problem which has remained open for more than $70$ years. It states that all graphs on at least three vertices are determined up to isomorphism by their deck \cite{k, MR0120127}. To know more about this problem, see \cite{MR1161466, MR0480189}. The deck of a graph $G$ is the multiset of graphs that is obtained from deleting one vertex in every possible way from the graph $G$. The members of a deck are referred to as \textit{cards}. A graph is  \textit{reconstructible}, if it is determined up to isomorphism by its deck. A class of graphs is reconstructible, if every member of the class is reconstructible. Any property of a graph, i.e. graph invariant, which is determined by the deck of graph, is said to be reconstructible. In this paper, the asymmetry of a graph, i.e. a subgraph which does not repeat, is employed to  reconstruct it form its deck.\\
The graph  reconstruction  conjecture was proposed by Ulam \cite{ MR0120127} in 1960. Kelly \cite{MR0087949}, in his Ph.D thesis, showed that regular graphs, Eulerian graphs, disconnected graphs and trees are reconstructible.  Graphs in which no two cycle have a common edge \cite{MR0252255}, graphs in which all  cycles pass through a common vertex \cite{MR0485581}, outer planner graphs \cite{MR0345873}, separable graphs without end vertices \cite{MR0262098}, maximal planar graphs \cite{MR615314}, critical blocks and graphs with some specific degrees sequence \cite{Nash} are some well known classes of reconstructible graphs. Another approach to this problem is attempt to find reconstructible graph invariants. The fundamental result in this area is  Kelly's Lemma \cite{MR0087949}. It states that the number of occurrence of any proper subgraph of a graph is reconstructible. Many results in graph reconstruction are based on this lemma. The generalization of this lemma for spanning subgraphs \cite{MR635955} asserts  that the number of any disconnected spanning subgraph is reconstructible. Some known reconstructible graph invariants are characteristic polynomial \cite{MR538033}, chromatic polynomial \cite{MR538033} and  planarity \cite{MR2344137}. Oliveira and Thatte \cite {MR3527999, MR2180800} had a new approach to this problem by studying the matrix of covering numbers of graphs by sequences of subgraphs and proposing a bound for the rank of this matrix.

In this paper, a general framework for the graph reconstruction problem is proposed. Here, we do not prove or reject the Reconstruction Conjecture. But, a general framework is proposed to explain the reconstruction of graphs. The reconstruction of small graphs have been shown by computer, while there was not any evidence that why small graphs are reconstructible. Here, we have explained  why small graphs are reconstructible.

We will show that any non-regular graph  has a proper induced subgraph which is  unique due to its structure or the way of its connection to  the rest of the graph (Lemma \ref{non-vertex}). The former subgraph is called an anchor and the latter is called a connectional anchor, if it is distinguishable in the deck.
%
%First, we show that how graph anchor helps us in graph reconstruction and introduces large families of reconstructible graphs.
%Then, we show that graph anchor enables us to have a general framework for the graph reconstruction problem.

In this paper, it is shown that a graph $G$  is reconstructible, if it has either an  orbit $O$ (with at least three vertices) which makes $G \backslash O$ an anchor or two vertices $\lbrace v,w  \rbrace$ which make $G \backslash \lbrace v,w\rbrace $
 an anchor under the conditions which will be mentioned. Also, we discuss the sufficiency of this result for any non-regular graph with at least three vertices.

 The above simple result is sufficient to prove that  many significant families of graphs are reconstructible, such as almost every graph and small graphs. %It may be wondering how this simple statement may be sufficient for every non-regular graph.
 What makes the above result, nearly, comprehensive, is the idea of  anchor extension. Accordingly, one approach to check the conjecture of graph reconstruction is to verify the above simple statement for all  non-regular graphs.%: Any graph has either an $(n-2)$-vertex (connectional) anchor  with the condition, which will be mentioned, or an orbit with at least three vertices whose removal leaves an anchor (or connectional anchor).

%By extending an anchor of a graph, we, usually, reach to an anchor whose removal leaves two vertices or an orbit with at least three vertices. If this event occurs always (Assumption I),  then every graph is reconstructible, provided that the mentioned condition for  an anchor whose removal leaves two vertices hold true (Assumption II). So far, no counterexample  has  been  found for these two assumptions. Accepting these two assumptions, the reconstruction conjecture holds true. In return, finding any counterexample for any of these assumptions may help us to find counterexample for the conjecture.

After preliminary definitions,  anchor and shadow are defined in  Section \ref{as} and  some results about anchor are given to have more intuition about it. In Section \ref{re}, the application of graph anchor for the reconstruction of graphs is described.  The  concepts of anchor extension and maximal anchor  are employed to draw a general framework for the graph reconstruction problem in Section \ref{maxi}.  In Section \ref{anchor_free}, we deal with  anchor-free graphs and define the concept of connectional anchor to explain the reconstruction of anchor-free graphs.  The reconstruction of graphs with ($n-2$)-vertex anchor is discussed in Section \ref{n-2}. The paper is concluded in Section \ref{conc}.
\section{Definitions and Notations}
 In this paper, any graph is simple and any subgraph is  vertex induced subgraph, otherwise, it is mentioned. The vertices and edges of a graph $G$ are denoted by $V(G)$ and $E(G)$, respectively. The order of a graph is the number of its vertices. The degree and the neighbors of a vertex $v \in V(G)$  are denoted by $d(v)$ and $N(v)$, respectively. If a vertex $v$ is deleted from a graph $G$, we denote it by $G \backslash v$. The multiset of graphs that is obtained from deleting one vertex in every possible way from the graph $G$, is the deck of $G$. \\
A graph $ G $ is called \textit{asymmetric}, if it has no nontrivial symmetries, i.e. $\rm{Aut }(G)=I$. Two vertices $u$ and $v$ of a graph $G$ are called \textit{similar}, if there is an automorphism of $ G$ which maps $u$ into $v$. %Dissimilar vertices whose removal leaves isomorphic subgraphs are called \textit{pseudo-similar}\cite{MR1470790, MR615005}.
 Similarity is, obviously, an equivalence relation. Thus, the similar vertices are in classes which are called \textit{orbits}. We call the subset $S \subset V(G)$  invariant under $\rm{Aut}(G)$, if $\lbrace \theta(S) \vert  \theta \in \rm{Aut(G)} \rbrace=S $. If $H$ is a subgraph of a graph $G$, $G \backslash V(H)$ is the induced subgraph on $V(G)-V(H)$.
 The \textit{cover} of a subgraph $H$ in $G$ is the minimal induced subgraph of $G$ which includes all possible induced subgraphs isomorphic to $H$ in $G$.\\
 In the probability space of graphs on $n$ labeled vertices in which the edges are chosen independently, with probability $p=1/2$, we say that almost every graph $G$ has a property $Q$ if the probability that $G$ has $Q$ tends to 1 as $n \rightarrow \infty$.\\%\cite{MR1864966}
For graphs $F$ and $G$, we use $G \choose F$ to denote the number of vertex induced subgraph in $G$ which are isomorphic to $F$. For instance, if $F = K_{2}$, then $G \choose F$ is the number of edges of $G$.
It is worthy to review Kelly's Lemma, here.
\begin{lemma}(Kelly's Lemma)\\

For graphs $H$ and $G$,
\begin{itemize}
\item If $\vert V(H) \vert <\vert V(G) \vert $, then $G \choose H $ is reconstructible. \rm{\cite{MR0087949}}
\item If $H$ is disconnected and $V(H)=V(G)$, then $G \choose H $ is reconstructible. (Kelly's Lemma generalization for spanning disconnected subgraphs) \rm{\cite{MR538033}}
\end{itemize}
\end{lemma}
\section{Graph Anchor and Shadow}\label{as}
Symmetry occurs when some patterns repeat. In contrast, asymmetry arises when there exists a pattern which does not repeat and is unique. In this way, a unique subgraph $H$ of a graph $G$, which does not repeat, is a representative of graph asymmetry. We define such subgraph as an anchor of graph $G$.\\
Unique subgraph is a known key concept for the graph reconstruction problem. Bollob\'{a}s \cite{MR1037416} has employed graphs in which all ($n-2$) and ($n-3$)-vertex subgraphs are unique to show that almost every graph is reconstructible by three cards. Muller \cite{MR0441789} has shown almost every graph is reconstructible using graphs in which all (${{n}}/ {2}$)-vertex subgraphs are unique. Unique subgraphs also have been used by Chinn \cite{MR0280395} and Zhu \cite{MR1470804} to introduce some families of reconstructible graphs. Ramachandran \cite{MR716448} has employed the idea of unique subgraph for the digraph $N$-reconstruction.

\begin{definition}
A proper induced subgraph $H$ of a graph $G$ is an \textit{anchor}, if it occurs exactly once in graph $G$. In other words, $G \choose H$ $ =1 $. The induced subgraph on $V(G)-V(H)$ is the \textit{residue} of $H$. We denote it by $G \backslash V(H)$.

 The \textit{anchor number} of a graph $G$ is the minimum order of the  anchors of $G$. We denote it by \rm{Anch}(G).

\end{definition}

Unique subgraph, also, has been defined by Entringer and Erd\H{o}s \cite{MR0317990} and used by Harary and Schewenk \cite{MR0325447}. They have used the concept of unique subgraph as the spanning subgraphs which are unique. Here, in contrast, we deal with proper vertex induced subgraph.

An anchor is a unique subgraph and, therefore, is distinguishable in any card containing it. Therefore, an anchor in a graph, similar to  a real anchor which fixes a boat, fixes a part of some cards and makes it possible to compare them.\\

\begin{figure}[ht]
\centerline{\includegraphics[width=8.5cm]{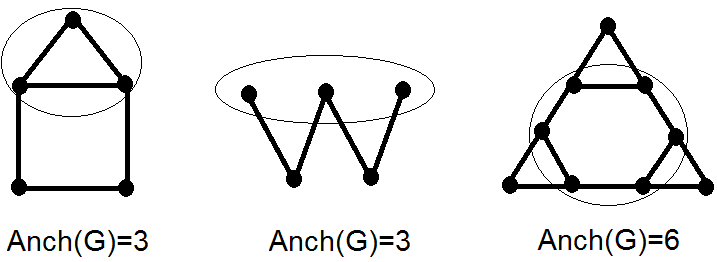}}
\caption{\label{exam}\small Anchor in some graphs }
\end{figure}

For illustration, some graphs with their anchors and  anchor numbers are shown in  Figure \ref{exam}.

Clearly, a graph may have some different anchors. In graphs which do not have any unique proper subgraph, we take the anchor number $n$, for example vertex-transitive graphs. The possible values for anchor number of a simple graph are $ 2, 3,.., n-1$ and $n$.  \\

\begin{table}
\begin{tabular}{c  c c c c c c  c c |c}
 $$ & $k$= 2 & $k$=  3 &$k$=4 &$k$=5&$k$=6 &$k$=7&$k$=8 &$k$=9 & Total\\
\hline
$n$=5 &2& 17& 12& 3& & & & &34\\
$n$=6 & 2&    48&    86 &    8 &   12& & &  &156\\
$ n$=7  & 2&   212   &654 &  146   & 26 &    4&&  &1044\\
 $n$=8 & 2 &       1062    &    7786     &   3082 &        373  &        10   &       31 & & 12346\\
$ n$=9 &   2 &       7266 &     139850&      121609&        5697 &        162 &         67 &         15&274668\\
 \end{tabular}
\caption{ \label{a} \small The number of $n$-vertex graphs with anchor number $k$}
\end{table}
To have an intuition about anchor number, the numbers of $n$-vertex graphs with anchor number $k$ are shown in Table \ref{a}. According to Table \ref{a}, for $ n$-vertex graphs where $n=7$, we have  \rm{Anch}$(G)<n-1$ for more than $97 \%$  of graphs.  For graphs on $n=8$ vertices,  more than 99.6 percent of graphs  posses an anchor with at most  ($ n-2$) vertices.
As we see in Table \ref{a} , 99.9 percent of $9$-vertex graphs have an anchor with at most 6 vertices.

Before dealing with graph reconstruction problem, we show some results about graph anchor.

\begin{proposition}
The subgraph $ H $  is an anchor of a graph $G$ if and only if $ \overline{H} $ is an anchor of $\overline{G} $. Therefore, the anchor numbers of $ G $ and $ \overline{G} $ are equal.

\end{proposition}

\begin{lemma}\label{simi}
An anchor either completely includes or excludes all
the vertices in each orbit. In other words, an anchor is an induced subgraph on the union of some orbits.
\end{lemma}
\textit{Proof}: If some vertices of one orbit are in an anchor while other vertices are not, the action of the automorphism group  makes another copy of the anchor, which contradicts the anchor definition. Therefore, an anchor either includes all vertices in an orbit $O$  or excludes any vertex in $O$. $ \diamond $

\begin{lemma}\label{}
In a graph $G$ with anchor number $\rm{Anch}(G)$, any possible induced subgraph $H$ of order $k$, where $ k<\rm{Anch}(G)$, repeats.
\end{lemma}
\textit{Proof}: Otherwise, $H$ is an anchor of $G$. Thus,   $\rm{Anch}(G) \leq k$. It contradicts to the assumption that   $ k<\rm{Anch}(G)$.$ \diamond$\\

If a graph is anchor-free, then any proper subgraph with any order is repeated. Under this condition,  a nearly symmetric structure is expected for such graphs. Due to Table \ref{a},  a few fraction of small graphs  are anchor-free. In a graph with \rm{Anch}($G$)=$n-1$, any subgraph with at most ($n-2$) vertices should be repeated. Accordingly, such graphs, similar to anchor-free graphs, are very rare.

\begin{lemma}\label{cover}
The cover of any arbitrary subgraph $H$ of $ G $, i.e. the minimal subgraph of $G$ which includes all possible copies of $H$ in $G$, includes at least $\rm{ Anch}(G) $ vertices.
\end{lemma}
\textit{Proof}: Let $H_c$ be the  minimal subgraph of $G$ which includes all possible copies of $H$, i.e. cover of $H$ in $G$.   Clearly, $H_c$ is an anchor of $G$, if it is a proper subgraph of $G$. Because, if there is another copy of $H_c$ in $G$, then there is a vertex out of $H_c$ which belongs to a copy of $H$. But, it contradicts to the assumption that $H_c$  includes all possible subgraphs isomorphic to $H$. Thus, $H_c$ is an anchor of $G$. Since, $Anch(G)$ is the minimum order of  anchors of $G$, we have $Anch(G) \leq \vert V(H_c) \vert$.$\diamond$
%Otherwise, the cover of $H$ in $ G $  is an anchor with smaller order, contradicting the definition of anchor number.
 %,  more symmetry is expected when the anchor number of a graph increases.
\begin{corollary}\label{cover_anchor_free}
Any possible subgraph of an anchor free graph covers all vertices of the graph.
\end{corollary}
\textit{Proof}:  If $G$ is an anchor free graph, then $Anch(G)=n$. Thus, according to Lemma \ref{cover}, for any arbitrary subgraph $H$ of $G$, we have $n \leq \vert V(H_c) \vert$. It means that the copies of $H$ covers all vertices of $G$. $ \diamond$\\
Here, we have a strong result about anchor-free graph. Any possible subgraph of an anchor free graph covers all vertices of the graph. Accordingly, a complete symmetric structure would be expected for anchor free graphs. The structure of small anchor free graphs confirms this expectation. Vertex-transitive graphs are an important family of anchor-free graphs. Anchor free graphs with at most 9 vertices are studied in Section \ref{anchor_free}.
%\begin{lemma}
%In a graph with no anchor, each vertex is a maximum degree vertex or adjacent to  a maximum degree vertex.\end{lemma}
%\textit{Proof}: Otherwise, the induced subgraph on vertices with maximum degree and their adjacent vertices is an anchor. Because, within the subgraphs with the same number of vertices it has maximum possible number of edges. $ \diamond $
%

\begin{figure}[ht]
\centerline{\includegraphics[width=3cm]{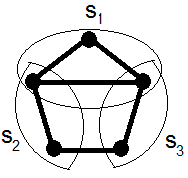}}
\caption{\label{shad}\small The shadows $s_1$ and $s_2$ are non-isomorphic. The shadows $s_2$ and $s_3$ are isomorphic. }
\end{figure}

\begin{definition}\label{}
A\textit{ shadow on a graph} is a subset of its vertices. Given a graph $G$, shadows $s_1$ and $s_2$ are \textit{isomorphic}, if there is an automorphism of $G$ that maps $s_1$ to $s_2$.
\end{definition}
In Figure \ref{shad}, three shadows $s_1$,$s_2$ and $s_3$ are shown in a graph. The shadows $s_1$ and $s_2$ are non-isomorphic, while $s_2$ and $s_3$ are isomorphic.
\begin{definition}
Let $H$ be a subgraph of a graph $G$ and $v \in V(G)-V(H)$. The neighbors of vertex $v$ in subgraph $H$, i.e. $N_{v} \cap H$, is the \textit{shadow of vertex $v$ on the subgraph $H$}.
We denote it by $s_{v,H}$ (or $s_v$ for simplicity, if $H$ is a specified subgraph).%$N_{v,H}$.
\end{definition}
So far, the concepts of shadow and isomorphic shadows were defined. Now, we define new concepts which relates to a set of shadows.

\begin{definition}
Let $S$ be a multiset of shadows on a graph $G$. We define a pair of shadows $s_i$ and $s_j$ are similar in $S$, if there exists $\theta \in Aut(G)$ such that $\theta (S)=S$ and $\theta(s_i)=s_j $. We call the multiset $S$ is shadow-transitive, if  any pair of shadows in $ S$ are similar.
\end{definition}

%--------------------------------------
%\begin{definition}
%Let $H$ be an anchor of a graph $G$ and $v$ a vertex out of the anchor. We call $N(v) \cap V(H)$ as the \textit{shadow} of $v$ on $H$ and denote it by $s_v$. Let $S$ be the set of shadows that vertices out of the anchor $H$ make on it.  Let $A,B \subset S$. We define $A$ and $B$ are isomorphic shadows, if there is $\theta \in Aut(H)$ such that $A=\theta(B)$. We call $A \subset S$ is an anchor of $S$, if there is not any $A' \subset S$  which is isomorphic to $A$.
%\end{definition}
%For illustration, Fig. \ref{shadsam} shows a sample  set of five shadows $S=  \lbrace s_1,s_2,s_3,s_4,s_5\rbrace $ on a sample graph.  $\lbrace s_1 , s_4 \rbrace$  is a subset of this shadows set which is isomorphic to  $\lbrace s_2, s_5\rbrace$. $\lbrace s_2,s_3,s_5\rbrace$ is an anchor for this shadows set.
%\begin{figure}[ht]
%\centerline{\includegraphics[width=4cm]{f2}}
%\caption {\label{shadsam} \small A sample set of shadows }
%\end{figure}
%-----------------------------
\section{  Graph anchor as a tool  for reconstruction }\label{re}
How an anchor of a graph $G$ helps us in the reconstruction of $G$ from its deck? Let $H$ be an anchor of $G$. The subgraph $H$ is  unique and, consequently, is the same in any card containing a copy of $H$. Thus, there is a one to one mapping between the vertices of $H$ in any card containing $H$. An anchor provides a partial correspondence between the vertices of cards including $H$. This partial correspondence benefits us to  reconstruct a graph from its deck.
%Using an anchor of a graph the reconstruction of that graph reduces to a smaller form of the reconstruction. In \cite {farhadian2016reconstruction}, it is shown that how using an anchor of a graph $G$, the reconstruction of $G$ reduces to the reconstruction of a shadow graph with smaller order.
\begin{lemma}\label{paie}
If a graph $G$ has an anchor, then it is distinguishable in any card containing it.
\end{lemma}
\textit{Proof}: Let $G$ be graph with $n$ vertices. According to Kelly's Lemma, the number of
occurrence of any induced subgraph $H$ with at most
($n-1$) vertices is reconstructible from the deck of $G$. Therefore, the anchors of a graph can be obtained from the deck. Since any anchor is a unique subgraph, it is distinguishable and the same in any card containing it.

\begin{definition}
Let $s$ be a shadow on a graph $H$. We call $s$ is fixed on $H$, if we have $\theta(s)=s$ for any $\theta \in Aut(H)$.
\end{definition}
In other words, a shadow on a graph is fixed, if it can not be moved by the automorphism group of the graph.
\begin{theorem}\label{fix}
Let $G$ be a graph with anchor $H=G \backslash \lbrace v, w \rbrace$.   If the shadow $s_{v,H}$ (or $s_{w,H}$)  is fixed  on $H$, then $G$ is reconstructible.
%If the neighbors of $v$ in $H$, that is $s_{v}$, is invariant under $\rm{Aut}(H)$, then $G$ is reconstructible.
\end{theorem}
\textit{Proof}: Due to Lemma \ref{paie}, anchor $H$ is distinguished in all cards containing it, that is in the cards  $G \backslash \lbrace v \rbrace$ and  $G \backslash \lbrace w \rbrace$. Thus,  vertex $v$ and its neighbors on $H$ are distinguished in $G \backslash \lbrace w \rbrace$.  Since, $s_{v}$ is fixed on $H$, there is just one realization for $s_v$ on $H$. Therefore, card $G \backslash \lbrace w \rbrace$ can be completed to $G$ by adding vertex $v$  adjacent to the set  $s_v$ in $H$. The existence of edge between $v$ and $w$ can be inferred from the number of edges which is reconstructible.$\diamond $

 In the previous theorem, it was assumed that the shadow of vertex $v$ on $H$ is fixed under the automorphism group of anchor $H$. In the next theorem, the shadows $s_v$ and $s_w$ are not fixed under the automorphism group of anchor. But, the distance between $v$ and $w$ specifies  the state of placing both of $s_v$ and $s_w$  together in $H$.

\begin{theorem} \label{distance_anchor}
Let $G$ be a graph with anchor $H=G \backslash \lbrace v, w \rbrace$. If the distance of $v$ and $w$ within $H$ specifies the state of placing both of the neighbors of $ v$ and $w$ together in $H$, then graph $G$  is reconstructible.
\end{theorem}
\textit{Proof}: The neighbors of $w$ and $v$ in $H$  are distinguishable in cards $G\backslash v$ and $G\backslash w$, respectively, up to isomorphism. It is sufficient to know the position of them when they come together. According to the hypothesis, the distance of $v$ and $w$ within $H$ clarifies the relative position of the neighbors of $v$ and $w$ in $H$. Thus, it is sufficient to show the distance of $v$ and $w$ within $H$ is reconstructible.\\
The number of subgraphs containing $v$ and the number of subgraphs containing $w$ can  be obtained from the cards $G\backslash w$ and $G\backslash v$, respectively. In addition, the subgraphs which include  none of them are, exactly, the subgraphs of $H$ and, thus, their numbers is reconstructible. Therefore, the number of subgraphs which include both $v$ and $w$ are reconstructible. The smallest path which includes $v$ and $w$ indicates the distance of $v$ and $w$ in $H$ when $v$ and $w$ are not adjacent in $G$. If $v$ and $w$ are adjacent in $G$, we consider the smallest cycle including $v$ and $w$.  Thus, the distance of $v$ and $w$ within $H$ is reconstructible. $\diamond$
\begin{theorem}\label{}
Let $G$ be a graph with asymmetric anchor $H$ and at least two vertices in the residue. If there is no pair of vertices in $V(G)-V(H)$ with the same shadow on $H$, then $G$ is reconstructible.
\end{theorem}
\textit{Proof}:
Anchor $H$ is asymmetric. Thus, we have  $\vert \rm{Aut}(H)\vert =1$. It means  that all  shadows are fixed on $H$. If $\vert V(G)-V(H) \vert=2$, then graph is reconstructible due to Theorem \ref {fix}. Thus, suppose $\vert V(G)-V(H) \vert>2$.\\
 The subgraph $H$ is asymmetric and no pairs of vertices in $V(G)-V(H)$ have the same neighbors in $H$.  Thus, any vertex of residue is identifiable by its shadow on the anchor in the cards containing the anchor $H$. Thus, we can recognize adjacent and non-adjacent pair of vertices in residue. Therefore, we add the residue vertices  to anchor, using their shadows and make connection between the adjacent pairs.$\diamond$

\begin{theorem}\label{thm:twoA}
If the induced subgraph on the neighbors (or non-neighbors) of a vertex of a graph $G$  is an anchor, then graph $G$ is reconstructible.
\end{theorem}
\textit{Proof}: Let $v$ be a vertex of graph $G$ in which the induced subgraph on the neighbors of $v$, say $H$, is an anchor of $G$. If $d(v)=n-1$, then $v$ is an isolated vertex of $\overline{G}$ and graph $\overline{G}$ is reconstructible due to being disconnected. Consequently, $G$ is reconstructible. Now, suppose that $d(v)\leq (n-2)$. %Therefore, $\vert V(H) \vert \leq (n-2)$ and $H$ is present in at least two cards. The anchor $H$ is distinguishable in any card containing it according to Lemma \ref{paie}. We gather all such cards. If we omit $H$ from these cards, we have the deck of $G \backslash V(H)$.
Since all neighbors of $v$ belong to $V(H)$, $v$ is an isolated vertex in $G\backslash V(H)$. Let $H'=G\backslash V(H)$. The disjoint union of $H$ and  $H'$, that is $H \cup H'$, is a disconnected spanning subgraph. Thus,  the number of its occurrence is reconstructible according to the generalization of Kelly's Lemma. Since there is just one subgraph $H$ in $G$, there is just one spanning subgraph isomorphic to $H \cup H'$. If there is a subgraph $H"$ such that the number of occurrence $H \cup H"$ is not zero, then $H"$ is isomorphic to $H'$ or  one of its subgraphs. Therefore,  subgraph $H'$, i.e. $G \backslash V(H)$,  is reconstructible from the deck. Clearly, vertex $v$ is an isolated vertex in subgraph $H'$.  Within cards containing anchor $H$, we choose cards in which an isolated vertex of $H'$ is deleted from.  Within the isolated vertices of $H'$, there is at least one vertex which is adjacent to all vertices of $H$, that is $v$.  We know that for any card $G\backslash \lbrace x \rbrace$, the degree of vertex $x$ is reconstructible. Vertex $v$ is adjacent to all vertices of $H$. Thus, it has maximum degree in $G$ within all isolated vertices of $H'$.  In a card that  an isolated vertex of $H'$ with maximum degree in $G$ has been deleted, a vertex should be added adjacent to all vertices of $H$ and isolated from all vertices of $H' \backslash v$. Therefore,  $G$ is reconstructible. $\diamond$
%In a card that $x$ is deleted from, . But,
%In card $G-\lbrace x \rbrace$, we add $x$ without any edge to vertices of $G \backslash V(H)$.
%The number of graph edges, which is reconstructible from the deck, imposes at least one of them, i.e. vertex $v$, to be adjacent to all vertices of $H$. Thus graph $G$ is reconstructed  by add $v$ without any neighbor in $G \backslash V(H)$ and adjacent to all vertices of $H$.

The previous results are more applicable for graphs which are asymmetric or close to asymmetric graphs in structure. The following theorem is useful for graphs with non-trivial automorphism group.

We know that if a graph $G$  is vertex-transitive, all cards of the deck of $G$ are the same. In opposite, if all cards of a graph $G$ are the same, then $G$ is vertex-transitive and, also is reconstructible due to regularity. We have a similar result for  the set of  shadows. If the shadows set $S$ is  shadow-transitive, i.e. all shadows  in $S$ are similar, then  all one shadow deleted subset of $S$, i.e $S-\lbrace s_i\rbrace$, are isomorphic. Following lemma shows that the inverse is also true. That is, if all one shadow deleted subset of $S$ are isomorphic, then $S$ is shadow-transitive. And the shadow set $S$ can be reconstructed from the multiset of one shadow deleted subset of $S$, i.e. $\lbrace S- \lbrace s_i \rbrace \vert s_i \in S \rbrace $.

%Clearly, if $S$ is a set of shadws then   for any pair of shadows $s_i$ and $s_j$ in $S$, we have $S-\lbrace s_i \rbrace=\theta (S- \lbrace s_j  \rbrace)$. It means that the deletin of any shadow reult isomorphic
\begin{lemma}\label{shadow_trans}
Let $G$  be an arbitrary graph and $S$ be a set of shadows on $G$ with at least three shadows. If all one shadow deleted subset of $S$ are isomorphic, then $S$ is reconstructible from multiset  $\lbrace S- \lbrace s_i \rbrace \vert s_i \in S \rbrace $.
%If $S-\lbrace s_i \rbrace$ and $S-\lbrace s_j \rbrace$ are isomorphic subset of shadows  on $G$ for any pair of shadows $s_i, s_j \in S$, then $S$ is transitive and can be reconstructred from  multiset $\lbrace S- \lbrace s_i \rbrace \vert s_i \in S \rbrace $.
\end{lemma}
\textit{Proof}: Let $D_S$ be the multiset of all one shadow deleted subset of $S$, i.e $D_S=\lbrace S- \lbrace s_i \rbrace \vert s_i \in S \rbrace$. According to the assumption, for any pair of shadows $s$, $s'$ in $S$, $ S- \lbrace s \rbrace$  and $ S- \lbrace s' \rbrace$  are isomorphic. % It means $D_s$ includes just one distinct element which is repeated $\vert S\vert$ times.
 If there exists a pair of non-isomorphic shadows  $s, s' \in S$ on $G$, then $S -\lbrace s \rbrace $ can not be isomorphic to  $S -\lbrace s' \rbrace $. But, it contradicts to the above  assumption. Thus,  all shadows in $S$ are isomorphic in $G$.  Since all shadows in $S$ are isomorphic in $G$, for any pair $s$ and $s' $ in $S$ there is $\theta \in \rm{Aut}(G)$ such that $s'=\theta(s)$. We draw an arc from $s$ to $s'$ labeled with $\theta$. Thus, we have a labeled directed graph $F_S$  such that  all cards of its deck are the same. Thus, graph $F_S$ is a vertex-transitive graph. Thus, in a card $S- \lbrace s \rbrace $ the deleted shadow $s$ can be added as $\theta(s_0)$ for a shadow $s_0 \in S- \lbrace s \rbrace$.

\begin{theorem}\label{orbit}
Let $O$ be an orbit of a graph $G$ with at least three vertices. If $G \backslash O$ is an anchor, then $G$ is reconstructible.
\end{theorem}
\begin{figure}[ht]
\centerline{\includegraphics[width=6.7cm]{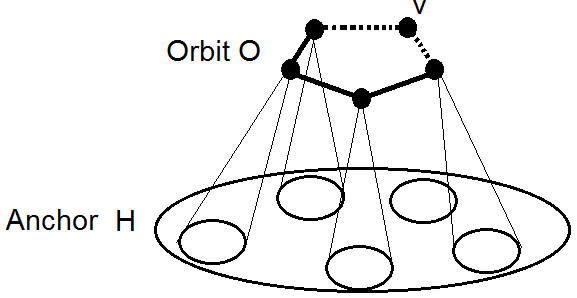}}
\caption{\label{Anchor_example}\small Graph $G$ with orbit $O$ and anchor $G\backslash O$ is reconstructible.}
\end{figure}
\textit{Proof}:
Let $H=G \backslash O$ be an anchor of graph $G$. According to Lemma \ref{paie}, anchor $H$ is distinguishable in any card containing it. We choose all cards containing anchor $H$. Since all vertices which are out of the anchor are similar, i.e. belong to the same orbit, all chosen cards are the same.  We distinguish anchor $H$ in all of them. If we omit anchor $H$ from these cards we have the deck of the induced subgraph on $G \backslash H$. Sine all chosen cards are the same, the deck of $G\backslash H$ includes a set of identical cards. Since, $G \backslash H$ is a regular graph,  vertex $v$ finds its neighbors in subgraph $G \backslash H - v$.  It is sufficient to find the neighbors of $v$ in the subgraph $H$, i.e. $s_{v,H}$. \\
Since, all cards containing the anchor $H$ are the same, all one shadow deleted subset on the anchor $H$ are the same, that is they  are isomorphic. Thus, according to Lemma \ref{shadow_trans}, the set of shadows on $H$ is reconstructible and vertex $v$  finds its neighbors on $H$. $\diamond$
\begin{remark}
Let $G$ be a graph with $n$ vertices. If all cards of the deck  $G$ are the same, we can surely conclude that all vertices of the graph are similar. %In contrast, if there are $k(<n)$ identical cards in the deck, we can not say their corresponding vertices are similar.
 If $G$ has an anchor $H$ with ($n-k$) vertices such that all $k$ cards containing the anchor $H$ are the same, we can surely conclude that $k$ vertices, which are out of the anchor, are similar due to Lemma \ref{shadow_trans}. Thus, graph $G$ is reconstructible due to Theorem \ref{orbit}.
\end{remark}

\subsection{Why small graphs are reconstructible?}
McKay \cite{MR1448235} has shown by computer that the reconstruction  conjecture is correct for graphs on at most 11 vertices,
but there was not any evidence why they are reconstructible.
Using the obtained results of the previous section, we can show the reconstruction of small graphs by reason. %As we see in the Table \ref {number_small_graph}, more than 99 percent of graphs on $n$ vertices where $n<10$ have an anchor with at most $n-2$ vertex.  We use their anchor and the results of previous section to prove they are reconstructible.

\begin{table}
\begin{tabular}{c   c| c  c }

 Graph index &  Reason  &  Graph index & Reason\\
\hline
G14& Anchor-free&		 G102& Theorem \ref{thm:twoA}\\
G30 & \rm{Anch}($G$)=$n-1$& G103& Theorem \ref{thm:twoA}   \\
 G31 & Theorem \ref{distance_anchor} &  G104& Theorem \ref{distance_anchor}\\
G34 & Theorem \ref{distance_anchor}&  G112& Theorem \ref{distance_anchor}     \\
 G35 & Theorem \ref{distance_anchor}& G113& Theorem \ref{distance_anchor}    \\
G36&Theorem \ref{thm:twoA}&  G114& Theorem \ref{distance_anchor} \\
G37&Theorem \ref{thm:twoA}&  G115& Theorem \ref{distance_anchor}\\
G78& \rm{Anch}($G$)=$n-1$&   G118& Theorem \ref{distance_anchor}\\
G79&Theorem \ref{distance_anchor}& G119&  Theorem \ref{thm:twoA}\\
G80&Theorem \ref{distance_anchor}&  G120&Theorem \ref{thm:twoA} \\
G81&Theorem \ref{thm:twoA}& G121   &Theorem \ref{thm:twoA}\\
G83&Anchor-free&  G122 &        Theorem \ref{distance_anchor}\\
G93&Theorem \ref{distance_anchor}&  G123&        Theorem \ref{thm:twoA} \\
G94&Theorem \ref{orbit}& G124&Theorem \ref{thm:twoA}\\
 G95& Theorem \ref{distance_anchor}&G125&Theorem \ref{distance_anchor}\\
 G96 &  Theorem \ref{distance_anchor} &G126& Theorem \ref{thm:twoA}\\
 G97   &Theorem \ref{distance_anchor} &G127&Theorem \ref{thm:twoA}\\
G98    &    Theorem \ref{distance_anchor}     & G128& Theorem \ref{distance_anchor}\\
 G99 &     Theorem \ref{distance_anchor}&   G129& Theorem \ref{distance_anchor}\\
G100& Theorem \ref{distance_anchor} &G130& Anchor-free\\
 \end{tabular}
\caption{ \label{reason} \small The reconstruction of graphs shown in Fig. \ref{6ver} are explained using the obtained results. Please note that for some graphs the theorem is applied for the complement of the graph.}
\end{table}

\begin{figure}[ht]
\centerline{\includegraphics[width=11.2 cm]{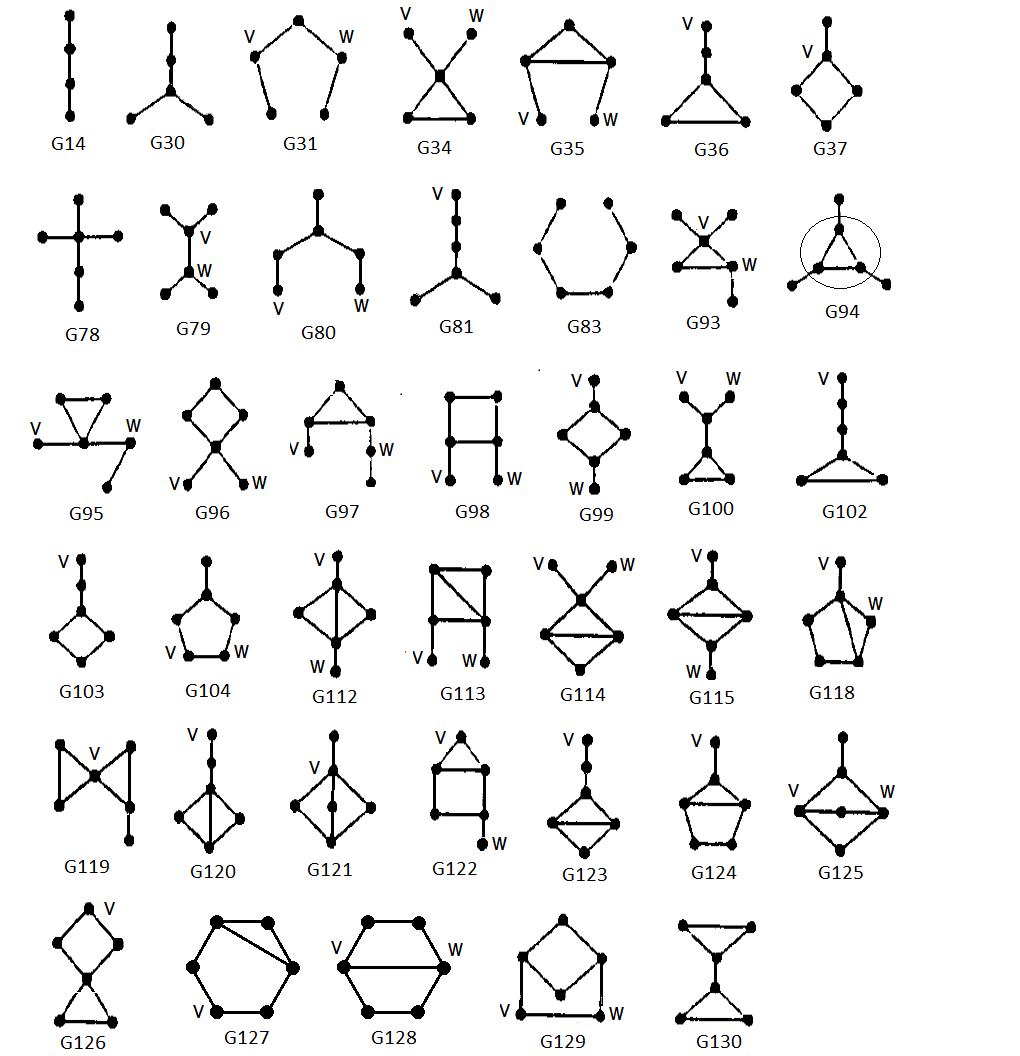}}
\caption{\label{6ver}\small  Non-regular connected graphs with at most 6 vertices and at most 15 edges which their complement is connected. The vertices which are out of the anchor are denoted by $v$ and $w$. }
\end{figure}
The theorems of Section \ref{as} are applied for graphs with at most 6 vertices. A list of graphs with at most 6 vertices are shown in Fig \ref{6ver}.  Since we know that disconnected and regular graphs are reconstructible, such graphs (and their complement) are omitted from this list. In addition, graph $G$ is reconstructible if and only if $G^{c}$ is reconstructible. Thus, only $n$-vertex graphs with at most $ {{{n} \choose {2}}/2 }$ edges are considered. The list of graphs and their assigned index are taken from \cite {MR1692656}. \\
 The anchor and theorems which prove these graphs are reconstructible, are shown in Figure \ref {6ver} and Table \ref{reason}, respectively. For each graph the vertices which are out of the anchor are assigned with labels $v$ or $w$. Unlabeled vertices belong to the anchor. The theorems which prove the reconstruction of  these graphs are listed in the Table \ref{reason}.
Please note that for some graphs, the theorems are applied for the complement of graph.
As the same way, we  can apply the mentioned results to show the reconstruction of graphs with 7 vertices or more. \\

 According to computation of anchor number for graphs with at most 9 vertices  in Table \ref {a} and due to symmetric structure, a few fraction of graphs are anchor-free  or with anchor number ($n-1$). The reconstruction of  small graphs with anchor number of $n$ and ($n-1$) is explained in Section \ref{anchor_free}

\section{Towards a general framework} \label{maxi}
\begin{figure}[ht]
\centerline{\includegraphics[width=3cm]{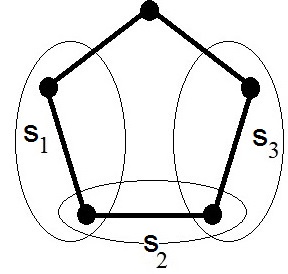}}
\caption{\label{neighbor_anchor}\small  $S=\lbrace s_{1}, s_{2}, s_{3}\rbrace $ is a set of shadows on $C_5$. $A=\lbrace s_{1}, s_{3}\rbrace $ is a subset of $S$.  There is  not any subset of $S$ which is isomorphic to $A $. Thus, the subset $A$ is an anchor for the shadow set $S$.}
\end{figure}

In the previous sections, some new results for the reconstruction of graphs were proposed. Now, we want to find a comprehensive approach for  graph reconstruction. In this way, the concepts of anchor extension and maximal anchor are defined in this section.
 As the main result of this section, we show that if an anchor is maximal, then its  residue  is anchor free. This fact is important. Because,  any arbitrary  anchor of any graph with anchor number $k$ where $k<(n-2)$, can be extended as long as we reach to an anchor free graph in the residue or just two vertices in the residue. Consequently, the reconstruction of any graph with $n$ vertices  reduces to the reconstruction of a graph with an anchor and anchor free residue or a graph with ($n-2$) vertices. As we will see in the next section, the anchor free graphs are very rare and have nearly symmetric structure. Vertex-transitive graphs are the majority of this family.\\
Before stating the main result, we define the required new concepts. \\
We want to  generalize the concept of anchor to the set of shadows. At first,  we define isomorphism for the set of  shadows.
\begin{definition}\label{}
Let $S=\lbrace s_{1}, s_2,\cdots, s_k \rbrace$ be a finite multiset of shadows on a graph $G$ and suppose that $A$ and $B$ are two arbitrary subsets of $S$.   Two shadow sets $A$ and $ B$ are \textit{isomorphic}, if there is an automorphism of $G$ that maps the subset $A$ to the subset $B$.

%Two finite set of shadows $S=\lbrace s_{1}, s_2,\cdots, s_k \rbrace$ and $S'=\lbrace s'_{1}, s'_2,\cdots, s'_k\rbrace$ on a graph $G$ are \textit{isomorphic}, if there is an automorphism of $G$ that maps the set  $S$ to the set $S'$.
\end{definition}
Now, we generalize the concept of anchor to the set of shadows.
\begin{definition}\label{}
Let $S$ be a set of shadows on a graph $G$. A subset $A$ of $S$ is a \textit{shadow anchor}, if there is not any subset of $S$ which is isomorphic to A.\end{definition}
In Fig. \ref{neighbor_anchor}, a graph with a set of shadows are shown. $S=\lbrace s_{1}, s_{2}, s_{3}\rbrace $ is a set of shadows on $C_5$. There is  not any subset of $S$ which is isomorphic to  subset $A=\lbrace s_{1}, s_{3}\rbrace $. Thus, the subset $A$ is an anchor for the shadow set $S$.

An anchor is a subgraph which is distinguishable in any card containing it, due to Lemma \ref{paie}. An anchor  provides a correspondence between the vertices of the anchor in different cards. Therefore, it is desired to have a larger anchor and, consequently, more information about the correspondence of cards  vertices.
\begin{definition}
Let $H_0$ be an anchor of a graph $G$. If we can find an anchor $H_1$ such that $V(H_0) \subset V(H_1)$,  we call anchor $H_0$ is \textit{extended} to $H_1$.\\
 An anchor is\textit{ maximal}, if it can not be extended to a larger anchor.

\end{definition}

%
%\begin{lemma}\label{psudo-set}
%Let $H$ be an anchor of graph $G$. If  $A$, $B$ are two subsets of $V(G)-V(H)$ with non-isomorphic shadows set on $H$, then the induced subgraphs on $V(H) \cup A$  and $V(H) \cup B$ are non-isomorphic.
%\end{lemma}
%\textit{Proof}: Let  $H_A$ and $H_B$ be the induced subgraphs on $V(H) \cup A$  and $V(H) \cup B$, respectively. Suppose they are isomorphic and there is an isomorphism from $H_A$ to $H_B$. The image of set $A$ is not set $B$, due to non-isomorphic shadows of $A$ and $B$ on $H$. Let $A'$ be the image of $A$ in  $H_B$. Therefore, deletion of  $A'$ from  $H_B$  is the same  $H_A$ where $A$ is deleted from, that is $H$. It means that $H_B-A'$ and $H$ are isomorphic, contradicting to the anchor definition.$\diamond$\\

Now, we are ready to state the main result of this section.

\begin{theorem}\label{maximal_anchor}
If an anchor is maximal, then its residue is anchor free. Also, the set of shadows that the vertices of its  residue establish on the anchor is anchor free.
%If $H$ is a maximal anchor of graph $G$, subgraph $G \backslash V(H)$ and the shadows that the vertices $V(G)-V(H)$ produce on $H$ are anchor free.
\end{theorem}
The above theorem states that if an anchor is maximal, we can not find any anchor in the residue. In addition, the shadows set that the residue vertices make on  a maximal anchor is anchor-free. To prove this theorem, we need the following lemma.

\begin {lemma}\label{extend}
 Let $H$ be an arbitrary anchor  of a graph $G$.
\item [a)] If $H'$ is an anchor of  $G \backslash V(H)$, then the induced subgraph on $V(H) \cup V(H')$ is an anchor of $G$.
\item [b)] Let $S_H$ be the multiset of shadows that the vertices of $G \backslash V(H)$ make on $H$. If there exists $A \subset V(G)-V(H)$ such that  the shadows of $A $ on $H$ is an anchor of $S_H$, then the induced subgraph on $V(H) \cup A$ is an anchor of $G$.

%Let  $A \subset V(G)-V(H)$. If  the induced subgraph on $A \subset V( G\backslash V(H)$ or its shadow on $H$ is  an anchor, then the  induced subgraph on $V(H) \cup A$ is an anchor of graph $G$
\end{lemma}
\textit{Proof}:
a) Let $H'$ be an anchor of $G \backslash V(H)$. We want to show that the induced subgraph on $V(H) \cup V(H')$ is also an anchor of $G$.   Suppose that the induced subgraph on $V(H) \cup V(H')$ is not an anchor of $G$, i.e. there is another copy of it. Please note that another copy of it should include $V(H)$. Because, subgraph $H$ is a unique subgraph.  Thus, suppose that there is $B \subset V(G)-V(H)$ that the induced subgraph on $V(H) \cup V(H') $ and $V(H) \cup B$ are isomorphic by mapping $\phi : V(H) \cup V(H') \to V(H) \cup B$.
Since, $H$ is an anchor and unique, $\phi(H)=H$ and $\phi(V(H'))=B$.  Therefore, the induced subgraph on $V(H')$, i.e. $H'$,  is isomorphic to the induced subgraph on $B$, contradicting to that $H'$  is an anchor of $G\backslash V(H)$.\\
b)  Let $A$ be a subset of $V(G)-V(H)$ whose shadows on $H$  is an anchor of $S_H$. Assume that the induced subgraph on $V(H) \cup A$ is not an anchor of $G$. Thus, there exists $B \subset V(G)-V(H)$ such  that the induced subgraph on $V(H) \cup A$ and $V(H) \cup B$ are isomorphic by mapping $\phi : V(H) \cup A \to V(H) \cup B$. Since, $H$ is an anchor and unique, $\phi(H)=H$ and $\phi(A)=B$. Thus, the restriction of $\phi$ to $V(H)$ is an automorphism of $H$. Therefore, the  shadow of $A$ on $H$ is mapped to the shadow of $B$ on $H$ by $\phi \in Aut(H)$ . It means that the shadow of $A$ is isomorphic to the shadow of $B$. But, it  contradicts to the assumption that the shadow of $A$ is an anchor of $S_H$.$\diamond$

The above lemma shows that an anchor can be extended  as long as there exists an anchor in the residue or a shadow anchor in the set of shadows on the anchor. Now, we are ready to prove Theorem \ref{maximal_anchor}

\textit{Proof of Theorem \ref{maximal_anchor}}:
Let $H$ be a maximal anchor of a graph $G$. Assume that  the residue has an anchor, say $H'$. According to Lemma \ref{extend}, the induced subgraph on $V(H) \cup V(H')$ is an anchor of $G$, contradicting to maximality of $H$. Therefore, the residue is an anchor free graph.\\
 Now, assume  $A \subset V(G)-V(H)$ produces a shadow anchor on $H$. Thus, the induced subgraph on $V(H) \cup A$ is an anchor of $G$ due to Lemma \ref{extend}. It contradicts to  maximality of $H$. Therefore, the shadows that the residue vertices produce on a maximal anchor is anchor-free.$\diamond$\\
Theorem \ref{maximal_anchor} has useful results. According to this theorem, any anchor of a graph can be extended until the residue and its shadows set does not have any anchor. Thus, any $n$-vertex graph with arbitrary anchor of order $k$ where $k<(n-2)$ can be extended as long as we reach to a large anchor of order ($n-2$) or  an anchor free residue.
\begin{remark}
 Theorem \ref{maximal_anchor} reduces the reconstruction of a graph to the reconstruction of an anchor free structure or a graph with a sufficiently large anchor, such as an $(n-2)$-vertex anchor.  Therefore, this theorem provides a general framework to deal with the reconstruction problem.
\end{remark}
 In the following sections, we deal with anchor free graphs and graphs with an anchor of order ($n-2$).\\
%Anchor free graphs have at most symmetry in the classification of graphs according to anchor number. The vertex transitive graphs are an important part of the family of  anchor free graphs. If the residue is a vertex transitive graph and all cards containing the anchor are the same, then graph is reconstructible according to Theorem \ref{orbit}. Anchor free graphs are investigated in Section \ref{anchor_free}.\\
%Suppose that there are $k$ vertices in the residue and  the anchor number of the residue is ($k-1$). Thus, the anchor expansion remains just one vertex in residue. Since such situation is not useful for graph reconstruction, we do not extend the anchor under this condition. Thus, it is necessary to study $n$-vertex graphs with anchor number ($n-1$), too. Section \ref{anchor_free} deals with these two families of graphs. \\
%%For graphs without any maximal anchor with at most $n-3$ vertices, we extend an anchor until  just two vertices remain in residue.
%According to Theorem \ref{maximal_anchor}, reconstruction of an arbitrary graph reduces to the reconstruction of an anchor free structure which will be studied in the Section \ref{anchor_free} or a graph with an ($n-2$)-vertex anchor. The reconstruction of graphs with an ($n-2$)-vertex anchor is investigated in the next section.

\section{Anchor-free graphs}\label{anchor_free}
%Anchor free graphs, i.e. graphs with anchor number $n$, have at most symmetry in the classification of graphs by anchor number.
We saw in Section 3 that  if $G$ is an anchor-free graph, then any arbitrary subgraph $H$ repeats. Moreover, the copies of $H$  cover all  vertices of the graph $G$. Thus, it would be expected that anchor free graphs, i.e. graphs with anchor number $n$, have most symmetric structure in the classification of graphs by the anchor number. The structure of small anchor free graphs confirm this expectation.\\
We saw in the previous section that the extension of anchor results in an anchor free graph in the residue. An anchor free graph is either formal, i.e. a vertex-transitive graph, or informal which will be discussed in this section. If residue and its shadows on the anchor are formal, in fact residue is an orbit of the graph and graph is reconstructible using Theorem \ref{orbit}.

 The main part of family  of anchor-free graphs are vertex-transitive graphs which are reconstructible due to regularity. Although, anchor free graphs which are not vertex-transitive are very rare, but they exist. Therefore, we should investigate anchor-free  graphs which are not vertex-transitive.
\begin{example}
Even paths and Cartesian product of an even path in a vertex transitive graph are examples of anchor-free graphs which are not vertex-transitive.
\end{example}
\begin{figure}[ht]
\centerline{\includegraphics[width=14cm]{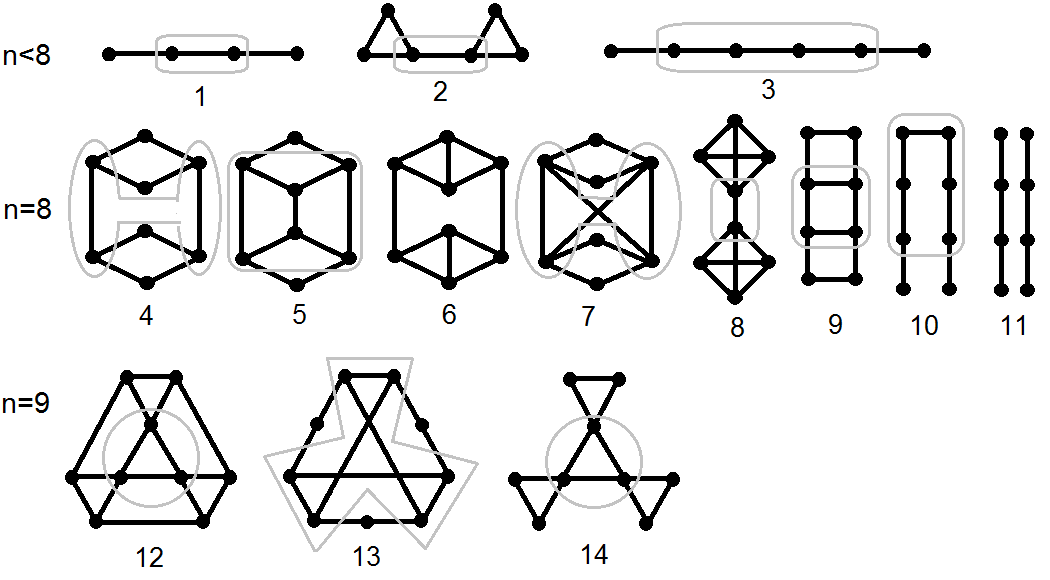}}
\caption{\label{anchor free}\small All anchor-free graphs which are not vertex transitive with at most 9 vertices}
\end{figure}

Small anchor-free graphs which are not vertex-transitive  with at most 9 vertices are shown in Figure \ref{anchor free}. Since \rm{Anch}($G$)=\rm{Anch}($G^c$), a graph or its complement is considered.

\begin{figure}[ht]
\centerline{\includegraphics[width=12cm]{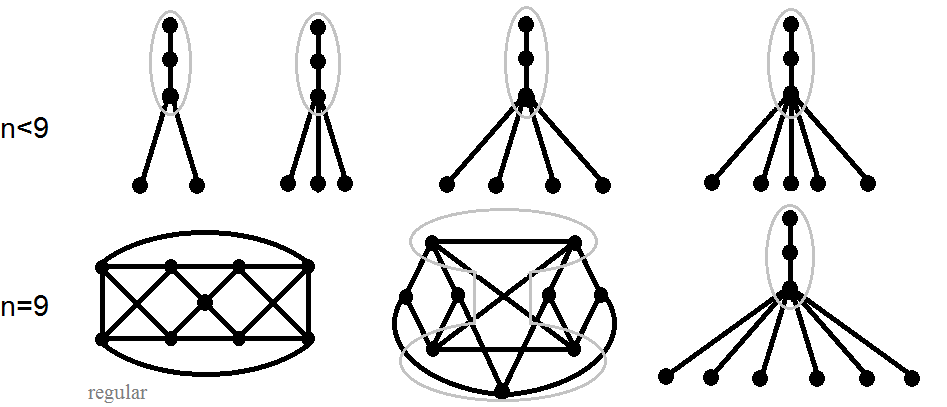}}
\caption{\label{n-1}\small  All $n$-vertex graphs with anchor number ($ n-1$) which are not disjoint union of an anchor free graph and an isolated vertex, with at most 9 vertices}
\end{figure}

If there are $k$ vertices in the residue and anchor number of residue is ($k-1$), then the extension of anchor remains just one vertex in the residue which is not useful for graph reconstruction. Therefore, the study of $n$-vertex graphs with anchor number ($n-1$) is required.
The structure of graphs with anchor number ($n-1$) is very close to anchor-free graphs.  The prevailing part of this family are disjoint union of an isolated vertex with an anchor free graph or complement of that, i.e. $K_1 \cup F$ or $\overline{K_1 \cup F}$  where $F$ is an anchor free graph. In this family,  a few number of graphs  have not such structure. The exceptional cases for at most 9 vertices, are shown in Figure \ref{n-1}. Since  \rm{Anch}($G$)=\rm{Anch}($G^c$), a graph or its complement is shown.

\begin{table}
\begin{tabular}{c  c   c c c}
  & \# & Anch(G)=$n$ or $n-1$ & Formal + Informal  \\
\hline
n=5 &34&15  (44 \%)&	12+2  \\
n=6 &156 &20  (13\%)& 16+4	\\
n=7 &1044&30    ( 3\%)&	28+2\\
n=8 &12346& 41 (0.03\%)&23+18\\
n=9 & 274668 &82 (0.03 \%)& 70 +  12	 \\
 \end{tabular}
\caption{\label{number_anchorfree} \small The number of graphs with anchor number of $n$ or ($n-1$)}
\end{table}
The family of graphs with anchor number of $n$ or ($n-1$) includes a small fraction of $n$-vertex graphs. The number of such graphs with at most 9 vertices are shown in Table \ref{number_anchorfree}. Formal anchor free graphs means vertex transitive graphs. And Formal graphs with anchor number ($n-1$) are disjoint union of an anchor free graph and an isolated vertex.\\

Here, we want to study the reconstruction of these two families of graphs. The prevailing part of anchor-free graphs are vertex-transitive which are reconstructible due to regularity. Therefore, it is sufficient to investigate non-vertex-transitive anchor free graphs for reconstruction.  Also, for graphs with anchor number of ($n-1$), the majority have an isolated vertex (themselves or their complement) which makes them reconstructible. \\
To explain the reconstruction of anchor free graphs or graphs with anchor number ($n-1$), the concept of \textit{connectional anchor }is defined in the following.
\subsection{Connectional Anchor}
%What makes an anchor as a useful tool for the  reconstruction of a graph from its deck, is the property of being distinguishable.
We know that  the unique structure of an anchor makes it distinguishable within all subgraphs, But is it possible to distinguish a subgraph  using something else?\\
The answer is yes. Suppose that in a graph $G$, there is  a vertex $v$ with degree $d$ such that no vertices has degree $d$ or $d-1$. So, it is just vertex $v$ whose degree is $d$ or $d-1$ in the deck. Thus,  vertex $v$ can be  distinguished in the deck by its degree.
In this example, it is not the structure of the vertex $v$ which makes it distinct. In fact, its connection to the rest of the graph makes it distinct.
We can use this idea to distinguish a subgraph. Suppose that there are two copies of subgraph $H$ in a graph $G$. For one of them, there are 10 edges between $H$ and $G \backslash V(H)$ and in any card of the deck which is containing this copy, there are at least 7 edges between $H$ and the rest of the graph. In opposite, for other copies of $H$, there are just 5 edges between $H$ and $G \backslash V(H)$. Therefore, the first copy is distinguishable in any card of the deck due to the number of edges between $H$ and $G \backslash V(H)$.\\
This example illustrates the idea of connectional anchor.
\begin{definition}\label{difr}
Let $H$ be  a subgraph of a graph $G$. We define \textit{the connection} of subgraph $H$ to the rest of graph by triple $(H, G\backslash V(H), S_H )$ where $S_H$ is the shadows set that $V(H)$ establish on  $G \backslash V(H)$\\
For two subgraphs $H$ and $H'$ in a graph $G$. We call the connection $H$ to the rest of graph is similar to the connection of  $H'$ to the rest of graph, if
\begin{itemize}
\item [I)]  Subgraph $H$ is isomorphic to subgraph $ H'$.
\item  [II)] Subgraph $G \backslash V(H)$ is isomorphic to subgraph $ G \backslash V(H')$.
\item [III)]Shadow set $S_H$ is isomorphic to shadow set $S_{H'}$.
\end{itemize}

% is different from the connection $H'$ to rest of graph, if either
%\begin{itemize}\label{difr}
%\item [I)]  $G \backslash V(H)$ is not  isomorphic to  $G \backslash V(H')$, or
%\item [II)]    The shadows set that $H$ establish on  $G \backslash V(H)$ is not isomorphic to the shadows set that $H'$ establishes on  $G \backslash V(H')$.
%\end{itemize}
\end{definition}

\begin{definition}
%Suppose that subgraph $H$ in a graph $G$ is repeated more than once.
We call  $H$ is a \textit{connectional anchor} of a graph $G$, if the connection of $H$ to the rest of graph is unique, i.e. there is not any subgraph  with isomorphic  connection. In addition, this distinction should be inferable from the deck.
\end{definition}

\begin{figure}[ht]
\centerline{\includegraphics[width=2.6cm]{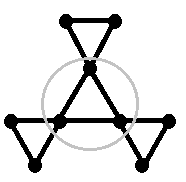}}
\caption{\label{co}\small  An anchor free graph with its connectional anchor  }
\end{figure}

\textit{Example 1}: In Fig.\ref{co}, an anchor-free graph is shown. This graph has four copies of $K_3$. Thus, $K_3$ is not an anchor. One of these $K_3$ is specified by a gray curve in Fig.\ref{co}. This copy of $K_3$ is different from others due to its connection to the rest of the graph. There are 6 edges between the specified $K_3$ and the rest of graph. If a vertex out of this subgraph is omitted, 5 edges will be remained. But, any other copies of $K_3$ has 2 edges in connection to the rest of graph.  Such copies have at most 2 edges in connection to the rest of graph in the deck. Thus, the specified $K_3$ is different from other copies due to the number of edges connecting it to the rest of graph and this distinction is distinguishable in the deck. Therefore, the specified $K_3$ is a connectional anchor and this graph is reconstructible using Theorem \ref{orbit}.$\square$

Similar to unique structure which makes a subgraph  distinguishable in the deck, the unique connection of a subgraph makes it distinguishable in the deck.
 The generalization of the Kelly's lemma states that  the number of spanning disconnected subgraph is reconstructible from the deck.
Since the  disjoint union of $H$ and $G \backslash V(H)$ is a disconnected spanning subgraph, the first and second condition of definition \ref{difr} can be checked by generalization of Kelly's lemma from the deck. %In any card containing $H$, $G \backslash V(H)$ is incomplete.\\%. However, the subgraph $H$ is a connectional anchor, provided that it is possible to distinguish $H$  in the deck using these  incomplete information about   $G \backslash V(H)$ in the card containing $H$. \\

For two different copies of  $H$ in a graph $G$, if two subgraphs  $G \backslash V(H)$ and $ G \backslash V(H')$ are isomorphic, then we should check the third condition of definition \ref{difr}. Since in any card containing $H$, $G \backslash V(H)$  is incomplete, i.e. one of its vertices is omitted, checking this condition is not easy. We have not a formal way to check this condition and depends on the structure of that graph.

Using connectional anchor, we are able to explain the  reconstruction of anchor free graphs. For small anchor-free graphs  given in  Fig.\ref{anchor free}, the connectional anchors which prove that they are reconstructible  are demonstrated by a gray curve.
Please note that any theorem and result based on graph anchor is, also,  correct  for connectional anchor. The theorems which prove the reconstruction of anchor free graphs of Fig \ref{anchor free}, are listed in Table \ref{reason_anchorfree}.
\begin{table}
\begin{tabular}{c   c| c  c }
 Graph index &  Reason  &  Graph index & Reason\\
\hline
1& Theorem \ref{distance_anchor} or \ref{Z2}&		 8& Theorem \ref{orbit}\\
2 & Theorem \ref{orbit}&	9 & Theorem \ref{orbit}   \\
 3 & Theorem \ref{distance_anchor} or \ref{Z2} &  10& Theorem \ref{distance_anchor} or \ref{Z2}\\
4 & Theorem \ref{orbit}&  11&  Disconnected  graph   \\
 5 & Theorem \ref{distance_anchor} or \ref{Z2}& 12& Theorem \ref{orbit}    \\
6&       Regular graph &  13& Theorem \ref{orbit} \\
7&Theorem \ref{orbit}&  14& Theorem \ref{orbit}\\
 \end{tabular}
\caption{ \label{reason_anchorfree} \small The reconstruction of  anchor free graphs shown in Fig \ref{anchor free} is explained using connectional anchor and the obtained results.}
\end{table}

The connectional anchors which explain the reconstruction of $n$-vertex graphs with anchor number ($n-1$) are also shown in Fig. \ref{n-1}. To prove the shown 5-vertex graph  in this figure is reconstructible, we can apply Theorem \ref{Z2} or \ref{distance_anchor}. For all other graphs shown in this figure, Theorem \ref{orbit} is sufficient.

\begin{lemma}\label{non-vertex}
Any non-regular graph has a proper induced subgraph which is  unique due to  its structure or the way of its connection to  the rest of the graph.
\end{lemma}
\textit{Proof:}  Let $G$ be an arbitrary non-regular graph and $H$  the induced subgraph on the vertices with maximum degree (or minimum degree). If there exists just one copy of $H$ in $G$, then $H$ is unique due to its structure. In opposite, if there are more than one copies of $H$, we show that $H$ is unique  due to its connection to the rest of graph $G$. We define function $f$ from subgraphs of $G$ to integer  numbers. We define $f(K)=2\vert E(K)\vert +\vert E(K,G\backslash V(K)) \vert$  where $K$ is a subgraph of $G$. We have
$ f(K)= \Sigma_{v \in V(K)} d(v)$. Within different copies of $H$ in graph  $G$, it is, exclusively, subgraph $H$ which posses maximum degrees. Thus, the value of $f(H)$ is maximum within all possible copies of $H$ in $G$. The value of $\vert E(H) \vert$ is the same for all  copies of $H$. Thus, $H$ is unique due to $\vert E(H,G\backslash V(H)) \vert$, i.e. the number of edges  between $H$ and $G \backslash V(H)$.$\diamond$
\begin{remark}
Please note that a subgraph which is unique due to its connection to the rest of the graph is not,  necessarily, a connectional anchor. This subgraph is a connectional anchor, if it is distinguishable in the deck. In contrast, a subgraph which is unique due to its structure is, always,  distinguishable in the deck.

\end{remark}

\begin{remark}
The definition of connectional anchor makes it possible to extend maximal anchor and have a larger anchor. The investigated examples of graphs show that an anchor can be extended until we reach to an orbit with at least three vertices or just two vertices  out of the anchor. But, we do not know whether  the anchor extension always can be continued to reach an orbit or two vertices out of the anchor.
\end{remark}

\section{Graphs with  ($n-2$)-vertex anchor}\label{n-2}
The reconstruction of graphs with an ($n-2$)-vertex anchor is important. Almost every graph has a ($n-2$)-vertex anchor. In addition, the anchor extension  results in an anchor free residue or two vertices out of the anchor. Here, we study the reconstruction of graphs in which the anchor extension remains  just two vertices out of the anchor, i.e. the reconstruction of $n$-vertex graphs with an ($n-2$)-vertex anchor.

Chinn \cite{MR0280395} has shown if there exists a vertex $v$ such that all ($n-2$)-vertex subgraphs of $G \backslash \lbrace v \rbrace$ are unique, graph $G$ is reconstructible. Zhu \cite{MR1470804} has improved this result by showing that at most three of ($n-2$)-vertex subgraphs of  $G \backslash \lbrace v \rbrace$  can be non-unique. Here, using Theorem \ref{fix} we have shown that one unique ($n-2$)-vertex subgraph which is asymmetric, is enough for $G$ to be reconstructible.

\begin{lemma}\label{almost}
Almost every graph has an asymmetric anchor with $(n-2)$ vertices.
\end{lemma}
\textit{Proof:} According to \cite{korshunov1985main, MR0441789, MR1037416}, for almost every $n$-vertex graph $G$,  all subgraphs with  $(n-3)$ vertices  are mutually non-isomorphic. In such graphs, any ($n-2$)-vertex subgraph ought to be unique and asymmetric.$\diamond$\\

\begin{theorem}\label{asymn-2}
Any $n$-vertex graph with an asymmetric anchor of order $(n-2)$ is reconstructible.
\end{theorem}
\textit{Proof}: Let $G$ be an $n$-vertex graph with asymmetric anchor $H$ of order ($n-2$). Since  $H$  is asymmetric, i.e. $ \vert Aut(H)\vert=1$, any shadow on graph  $H$ is fixed. Therefore, according to Theorem \ref{fix}, graph $G$ is reconstructible.

\begin{corollary}
Almost every graph is reconstructible.
\end{corollary}
\textit{Proof}: Using Lemma \ref{almost} and Theorem \ref{asymn-2}.$\diamond $
\\

In a graph $G$ with anchor $H=G \backslash \lbrace v,w \rbrace$, let $s_{v}$ and $s_{w}$ be the neighbors of $v$ and $w$ in $H$, respectively. $s_{v}$ and $s_{w}$ can be obtained from cards $G \backslash  w$ and $G \backslash  v  $, respectively. But, $s_{v}$ and $s_{w}$ can arbitrary move on $H$ by the action of \rm{Aut}($H$). Thus, it makes different possibilities for the state of placing both of them together and, consequently, two cards containing the anchor are not sufficient to reconstruct a graph, uniquely up to isomorphism.

\begin{figure}[ht]
\centerline{\includegraphics[width=10cm]{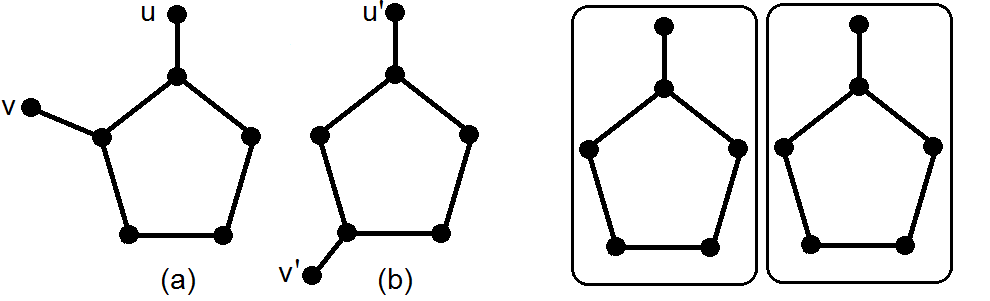}}
\caption{\label{relative}\small $s_{u,C_5}$ is isomorphic to $s_{u',C_5}$ and $s_{v,C_5}$ is isomorphic to $s_{v',C_5}$. But, the way of putting these two  shadows together on the anchor $C_5$ makes two non-isomorphic graphs. The set of cards containing the anchor $C_5$ is not sufficient to separate these two non-isomorphic graphs.  }
\end{figure}

For example, in Fig. \ref{relative}, two non-isomorphic graphs (a) and (b) are shown. The subgraph $5$-cycle is an anchor for both of them. The set of cards containing the anchor are the same and is shown in the right side. The set of cards containing the anchor is not sufficient to discriminate  two non-isomorphic graphs (a) and (b). Therefore, it is necessary to use other cards to find the relative  position of two shadows on the anchor. We do not know whether the other ($n-2$) cards are, always, sufficient to determine the position of two shadows on the anchor, together. The reconstruction conjecture claims that they are sufficient.

Theorem \ref{distance_anchor} uses the distance of vertices $v$ and $u$ to identify the state of placing both of $s_v$ and $s_u$ together, when the automorphism group of the anchor is not trivial. For example, in  Fig. \ref{relative} the distance between $v$ and $u$ in two graphs discriminates these two graphs. Thus, these two graphs are reconstructible due to this theorem.

In the following lemma, the order of automorphism group is restricted. Hence, the number of possibilities for putting two shadows together on the anchor is restricted.

\begin{lemma}\label{Z2}
Let $O$ be an orbit of a graph $G$ with two vertices such that $G \backslash O$ is an anchor and  the order of any element of  Aut($G \backslash O$) is at most three, then graph $G$ is reconstructible.
%Aut($G \backslash O$) $\in \lbrace I, Z_2, Z_3, D_3 \rbrace$
\end{lemma}
\textit{Proof}: Let $H= G \backslash O$ be an anchor of $G$ and  $O=\lbrace v_1 , v_2 \rbrace$ and $\lbrace s_{v_1} , s_{v_2} \rbrace$ be the shadows of $\lbrace v_1 , v_2 \rbrace$ on $H$. Since, $v_1$ and $v_2$ belong to the same orbit of $G$,  there exists $\theta \in \rm{Aut}(H)$ where $s_{v_1}= \theta(s_{v_2})$. Since according to assumption the order of $\theta$ is 1, 2 or 3, there are, exactly, two possibilities for $s_{v_1}$ and $s_{v_2}$: they are the same, or different. The degree sequence of the graph vertices which is reconstructible due to Kelly's Lemma, discriminates these two cases and chooses one of them. $\square$\\

%By anchor extension, we  reach to either  an anchor with ($n-2$) vertices or an anchor whose removal leaves an anchor free graphs. In this section, the former was studied, i.e. the reconstruction of graphs with an ($n-2$)-vertex anchor. Next section, we study the reconstruction of anchor free graphs.

\section{Conclusion}\label{conc}

In a graph, an anchor emerges when the graph symmetry is not enough to prevent it. Thus, except graphs with sufficient symmetry all graphs have an anchor. We saw that an anchor of a graph  can be an efficient tool to show that a graph is reconstructible.
 An anchor of a graph,  similar to a real anchor which fixes a boat, fixes a part of some cards to make it possible to compare them. Therefore, it benefits us for the  reconstruction of graphs.%enables us to infer how to add a vertex to the card which is deleted from.
\\
 It may be claimed that since  just cards containing the anchor are used for graph reconstruction, this method can not be efficient. But, this claim is not correct. Because, to find the anchors of a graph, we should consider subgraphs from all cards to find unique ones. Therefore, in this method we implicitly apply all cards to reconstruct a graph.\\
We saw that how graph anchor provides evidence for the reconstruction of small graphs. Small graphs, previously, had been shown to be reconstructible by computer.

In this paper, it is shown that  the anchor extension converts the reconstruction of any graph with anchor number $k$ where $k<(n-2)$ to the reconstruction of either an anchor free structure (or close to anchor free structure, i.e. with anchor number $n$ or $n-1$), or a graph with an $(n-2)$-vertex anchor.\\
 The reconstruction of graphs  with an $(n-2)$-vertex anchor was investigated and some new results were given. If a graph has an $(n-2)$-vertex anchor which is asymmetric, then  the graph is reconstructible. This fact results that almost every graph is reconstructible. It is sufficient to investigate graphs with an $(n-2)$-vertex anchor whose automorphism group is non-trivial. Generally, a graph with ($n-2$)-vertex anchor is reconstructible, if the relative position of two shadows on the anchor is reconstructible from the deck.\\

Also, the family of anchor free graphs and graphs with anchor number ($n-1$) were studied. The main part of anchor free graphs are vertex-transitive graphs which are reconstructible due to regularity.   To explain the reconstruction of non-regular anchor-free graphs, the concept of connectional anchor was introduced. We saw that the connectional anchor plays the same role of an anchor for the reconstruction of a graph from its deck. By applying the theorems and results obtained for anchor, to connectional anchor,  the reconstruction of anchor free graphs was explained. But, we have not shown that these method can be applied for every anchor-free graph. More study is suggested for the reconstruction of anchor free graphs. In addition, the definition of connectional anchor makes it possible to extend maximal anchor and have a larger anchor. The studied instances of graphs show that the anchor extension can be continued until we reach to an orbit with at least three vertices or just two vertices  out of the anchor. But, we do not know whether it is always possible.

In this paper, it is shown  that:\\
 If a graph $G$ has either an orbit with at least three vertices whose removal leaves an (connectional) anchor  or an ($n-2$)-vertex (connectional) anchor in which the relative position of  shadows is reconstructible, then $G$ is reconstructible.

The above simple statement is sufficient to explain the reconstruction of known reconstructible graphs, such as  small graphs and almost every graph. Thus, one approach to study the conjecture of graph reconstruction can be  the verification of  the following conjecture.
\begin{conjecture}
Any graph has either an $(n-2)$-vertex (connectional) anchor  in which the relative position of shadows on the anchor is reconstructible  or an orbit with at least three vertices whose removal leaves an anchor (or connectional anchor).
\end{conjecture}
 \bibliographystyle{plain}

%\bibliography{samplerefrences}

%
%
%\section*{Appendix \textit{I}: Why small graphs are reconstructible?}
%All graphs $G$ with at most 6 vertices such that $G$ and $G^c$ are not disconnected or regular, are shown in Fig. \ref{small}. For any graph, the anchor that makes it reconstructible, is shown by a gray closed curve.
%\begin{figure}[ht]
%\centerline{\includegraphics[width=10cm]{f3}}
%\caption{\label{small}\small All graphs with at most 6 vertices which are not disconnected or regular (also their complement) with the anchors  that prove they are reconstructible.  }
%\end{figure}
\end{document}